\documentclass[12pt,a4paper,reqno]{amsart} % 

\usepackage{microtype} % 自动优化字符间距和断行
\allowdisplaybreaks[4] % 允许多行公式跨页（1-4级）
\usepackage{slashed} % Slashed symbols e.g. to typeset Dirac operators
\usepackage{amsmath,amsfonts,amscd,amssymb,latexsym,mathtools}

\usepackage[left=2.9cm,right=2.9cm,top=2.66cm,bottom=2.66cm]{geometry}

\usepackage[all,matrix,arrow,tips,curve]{xy}
\usepackage[english]{babel}

\usepackage{mathrsfs}
\usepackage[percent]{overpic}

\usepackage{pgfplots}
\pgfplotsset{width=10cm,compat=1.9}
\usepackage{tikz}
\usetikzlibrary{shapes.geometric, calc}

\DeclareMathAlphabet{\pazocal}{OMS}{zplm}{m}{n}

\usepackage[nobysame,sorted-cites]{amsrefs} % alphabetic, msc-links,

\usepackage{enumitem}
\usepackage{comment}
\usepackage{soul}
\usepackage{marginnote}
\usepackage{bbm}
\usepackage{cancel,slashed}
\usepackage[normalem]{ulem}

\usepackage{footmisc}

\usepackage[colorlinks=true,allcolors=blue]{hyperref}
\usepackage[capitalise]{cleveref} %

\numberwithin{equation}{section}
\allowdisplaybreaks[1]

\theoremstyle{plain}% default 
\newtheorem{theorem}{Theorem}

\newtheorem{thmx}{Theorem}
\theoremstyle{definition}
\theoremstyle{definition}

\theoremstyle{remark} 
\newtheorem{remark}[theorem]{Remark}

\numberwithin{equation}{section}
\numberwithin{figure}{section}

\newcommand{\sphm}{\mathbf{S}^m} % m-sphere
\newcommand{\scal}{\operatorname{scal}} % scalar curvature
\newcommand{\supp}{\operatorname{supp}} % support
\makeatletter
\newsavebox{\@brx}
\newcommand{\llangle}[1][]{\savebox{\@brx}{\(\m@th{#1\langle}\)}%
\mathopen{\copy\@brx\kern-0.5\wd\@brx\usebox{\@brx}}}
\newcommand{\rrangle}[1][]{\savebox{\@brx}{\(\m@th{#1\rangle}\)}%
\mathclose{\copy\@brx\kern-0.5\wd\@brx\usebox{\@brx}}}
\makeatother
\makeatletter
\@namedef{subjclassname@2020}{\textup{2020} Mathematics Subject Classification}
\makeatother

\begin{document}

\title[Llarull's theorem on noncompact manifolds with boundary]{Llarull's theorem on noncompact manifolds with boundary}

\author[B. Liu]{Bo Liu}
\address{School of Mathematical Sciences\\ Key Laboratory of MEA(Ministry of Education) \\ \& Shanghai Key Laboratory of PMMP\\ East China Normal University\\ 200241, Shanghai\\ China}
\email{\href{mailto:bliu@math.ecnu.edu.cn}{bliu@math.ecnu.edu.cn}}
\urladdr{}
\thanks{}

\author[D. Liu]{Daoqiang Liu}
\address{Chern Institute of Mathematics \& LPMC, Nankai University, Tianjin 300071, China}
\email{\href{mailto:dqliu@nankai.edu.cn}{dqliu@nankai.edu.cn}}
\urladdr{\href{https://www.dqliu.cn}{www.dqliu.cn}}

\thanks{}

\date{\today}

\begin{abstract}
Recently, Zhang \cite{Zh20} and Li-Su-Wang-Zhang \cite{LSWZ24+} generalized Llarull's theorem to the noncompact complete spin manifold. In this paper,
we further extend their results to the noncompact manifold with compact boundary. 
\end{abstract}

\maketitle

%%%%%%%%%%%%%%%%%%%%%%%%%%%%%%%%%%%%%%%%%%%%%%%%%%%%%%%%%%%%%%%%%%%%%%%%%%%%%%
%%%%%%%%%%%%%%%%%%%%%%%%%%%%%%%%%%%%%%%%%%%%%%%%%%%%%%%%%%%%%%%%%%%%%%%%%%%%%%
%%%%%%%%%%%%%%%%%%%%%%%%%%%%%%%%%%%%%%%%%%%%%%%%%%%%%%%%%%%%%%%%%%%%%%%%%%%%%%

\section{Introduction} \label{sec:intro}
Llarull's sphere rigidity theorem \cite{Ll98} says that for a closed $m$-dimensional Riemannian spin manifold $(M,g_M)$ with scalar curvature satisfying $\scal_{g_M}\geq m(m-1)$, any smooth area decreasing map $\Phi:M\to \sphm$ of nonzero degree is a Riemannian isometry, where $\sphm$ is the standard unit $m$-sphere. 
For the noncompact extension of the Llarull's theorem, we summarize the recent works of Zhang \cite{Zh20} and Li-Su-Wang-Zhang \cite{LSWZ24+} as follows,
which completely answers a question of Gromov in an earlier version of \cite{Gr23}.
\begin{theorem}[{\cite{LSWZ24+},\cite{Zh20}}]\label{thm:NNSC_sphere}
Let $(M,g_M)$ be a complete noncompact Riemannian spin manifold.
Let $\Phi: M\to \sphm$ be a smooth area decreasing map which is locally constant at infinity and of nonzero degree.
Suppose that 
\begin{equation}
\scal_{g_M} \geq m(m-1)  \quad \text{on} \  \supp({\rm d}\Phi).
\end{equation} 
Then
\begin{equation}\label{eq:NNSC_sphere}
 \inf \scal_{g_M}< 0.
\end{equation}
\end{theorem}
Here a smooth map $\Phi:M\to \sphm$ is called area decreasing if for any two form $\alpha\in \Omega^2(\sphm)$, $\Phi^*\alpha\in \Omega^2(M)$ satisfies $|\Phi^*\alpha|\leq |\alpha|$.
We say that $\Phi$ is locally constant at infinity if it is locally constant outside a compact subset of $M$. Moreover, ${\rm d}\Phi$ is the differential of $\Phi$ and the support of ${\rm d}\Phi$ is defined to be $\supp({\rm d}\Phi):=\overline{\{p\in M\colon {\rm d}_p\Phi\neq 0\}}$.

Using the doubling procedure in \cite{deAl85} and \cite{GL80}, in this paper we extend 
\cref{thm:NNSC_sphere} to the noncompact complete manifold
with compact and mean-convex boundary.

\begin{thmx}\label{thmx:complete_long_neck_NNSC}
	Let $(M,g_M)$ be a complete noncompact $m$-dimensional Riemannian spin manifold with compact boundary.
	Let $\Phi:M\to \sphm$ be a smooth
	area decreasing map which is locally constant near the boundary $\partial M$ and at infinity and of nonzero degree.
	Suppose that
	\begin{enumerate}[label=\textup{(\arabic*)}]
		\item $\scal_{g_M} \geq m(m-1)$ on $\supp({\rm d}\Phi)$ and $\scal_{g_M}>0$ near $\partial M$; \label{complete_condition:scalar_comparison}
		\item $\partial M$ is mean-convex (i.e., the mean curvature $H_{\partial M} \geq 0$ with respect to the interior unit normal). \label{complete_condition:boundary_condition}
	\end{enumerate}
	Then %$\deg(\Phi)=0$.
	\begin{equation}%\label{eq:NNSC_sphere}
	\inf \scal_{g_M}< 0.
	\end{equation}
\end{thmx}
%\begin{remark}
	
The compact case of this result was studied in \cite{CZ24}*{Theorem~1.4} and \cite{Shi25}*{Theorem~1.1}.
%	Moreover, the even-dimensional case of \cref{thmx:complete_long_neck_NNSC} also follows from the techniques in \cite{Liu24+}. 
%	% A generalized result can be found in \cite{Liu24+}. 
%\end{remark}

%Using \cref{thmx:NNSC_Bends} and the previous similar doubling arguments, we have the following noncompact (necessarily not complete) long neck principle.

%, see also \cite{Liu24+} for an analogue.

%If $M$ is further not complete somewhere, we have the following version of the noncompact

% Recall that, for an exhaustion $K_1\subseteq K_2\subseteq\dots$ of a manifold $M$ by compact submanifolds, an \textit{end} of $M$ means a decreasing sequence $V_1\supseteq V_2\supseteq\cdots$ such that each $V_i$ is a connected component of $M\setminus K_i$.

\textbf{Notation.} Throughout this paper, we assume that all manifolds are smooth, connected, oriented, and of dimension no less than two. Moreover, the notion of completeness refers to metric completeness when $M$ has nonempty boundary.
% For the definitions of area decreasing maps and contraction constant $\|\Lambda^2 {\rm d}\Phi\|_{\infty}$ see \cite{Zh20}.

%%%%%%%%%%%%%%%%%%%%%%%%%%%%%%%%%%%%%%%%%%%%%%%%%%%%%%%%%%%%%%%%%%%%%%%%%%%%%%
%%%%%%%%%%%%%%%%%%%%%%%%%%%%%%%%%%%%%%%%%%%%%%%%%%%%%%%%%%%%%%%%%%%%%%%%%%%%%%
%%%%%%%%%%%%%%%%%%%%%%%%%%%%%%%%%%%%%%%%%%%%%%%%%%%%%%%%%%%%%%%%%%%%%%%%%%%%%%

\section{\texorpdfstring{Proof of \cref{thmx:complete_long_neck_NNSC}}{}}\label{sec:proof}

Since \(\partial M\) is compact and $\scal_{g_M}>0$ near $\partial M$, we can choose a compact collar neighborhood \(L = \partial M \times [0, \varepsilon_0]\) of $\partial M$, where $0<\varepsilon_0 \ll 1$, such that \(\scal_{g_M}>0\) on \(L\). By the compactness of \(L\), there exists \(\kappa > 0\) such that
\[
\scal_{g_M} \geq \kappa \quad \text{on } L.
\]
Since the argument in \cite{deAl85}*{Theorem 1.1} only depends on the local computations near boundary, it applies to noncompact manifolds with compact boundary.
By the trick in \cite{deAl85}*{Theorem 1.1}, it is enough to prove our theorem under the stronger assumption $H_{\partial M}>0$.

% To be specific, denote by $g_{\partial}$ the induced metric on $\partial M$, and assume that $\scal_{g_{\partial}} > 0$. Since $\partial M$ is compact, there exists $\kappa > 0$ such that
% \[
% \scal_{g_{\partial}} \geq \kappa \quad \text{on } \partial M.
% \]

% Following \cite{deAl85}*{Theorem 1.1}, fix a small number $\delta > 0$ (to be determined later) and choose a smooth function $h:\mathbf{R}\to \mathbf{R}$ such that
% \begin{itemize}
%     \item $h(t) = 1$ for $t \in [\delta, \infty)$;
%     \item $\| h - 1 \|_{C^2} < \delta$, i.e., $|h(t)-1|, |h'(t)|, |f''(t)| < \delta$ for all $t$;
%     \item $h'(t) < 0$ for $t < \delta$.
% \end{itemize}

% Consider a collar neighborhood $U\cong  \partial M \times [0, \varepsilon)$ of $\partial M$, where $0<\varepsilon\ll 1$, and equip it with a warped product metric of the form
% \[
% \widetilde{g}_M = h(t) g_{\partial} + dt^2.
% \]
% Then the scalar curvature of $\widetilde{g}_M$ is given by the standard formula for warped products:
% \[
% \scal_{\widetilde{g}_M} = \frac{\scal_{g_{\partial}}}{h} - (m-1) \frac{h''}{h} + \frac{(m-1)(4-m)}{4} \frac{(h')^2}{h^2},
% \]
% and the mean curvature of $(\partial M,\widetilde{g})$ satisfies
% \begin{equation}
%     \widetilde{H} = -\frac{m-1}{2} \frac{h'(0)}{h(0)} + H \geq -\frac{m-1}{2} \frac{h'(0)}{h(0)} >0.
% \end{equation}

Let $M_1:=M\backslash C$, where $C$ is a thin collar of $\partial M$.
Since $\partial M$ is compact,
we can choose $C$ and a small neighborhood $C'$ of $C$ in $M$ such that the mean curvature of $\partial M_1$
is still positive, the scalar curvature is positive on $C'$
and $\supp({\rm d}\Phi) \subset M\backslash
C'$. We consider the Riemannian product 
$M\times [0,1]$ and define
\[
D(M):=\{p\in M\times [0,1]: 
\mathrm{dist}(p, M_1\times \{1/2\})=\varepsilon \},
\]
where $0<\varepsilon\ll 1$. 
Then $D(M)$ is homeomorphc to the double of $M$.
Note that $D(M)$ can be divided into two parts: parallel part
$D_1:=(M\backslash C')\times [0,1]\cap D(M)$ and the bending part $D_2:=D(M)\backslash
D_1$.
By \cite{GL80}*{Theorem 5.7}, 
for $\varepsilon$ sufficiently small,
the metric on $D(M)$ can be smoothed such that the metric on $D_1$ is invariant
and the scalar curvature on $D_2$ is positive.
Note that the new metric on $D(M)$ is complete.	
%	
%%By contradiction, assume that $\deg(\Phi)\neq 0$. 
%Let us denote the double of $M$ by ${\rm d}M$. That is, ${\rm d}M:=M\cup_{\partial M} M^-$, where $M^-$ denotes the manifold $M$ with opposite orientation.
%Since $\scal_{g_M}>0$ on $M$ and $H_{\partial M}>0$ on $\partial M$, there exists a complete smooth Riemannian metric $g_{{\rm d}M}:=g_M\cup g_M$ on ${\rm dM}$ such that $\scal_{ g_{{\rm d}M} } > 0$ on ${\rm d}M$ and $g_{{\rm d}M}$ is actually a two-copies of $g_M$ outside a small enough neighborhood of $\partial M$ by \cite[Theorem~5.7]{GL80}(the arguments of Gromov-Lawson actually also apply to complete manifolds with compact strictly convex boundary).
%Note that ${\rm d}M$ is an open (i.e., noncompact, without boundary) 
The manifold $D(M)$ is equipped with a spin structure induced from the spin structure of $M$.

Fix $\star \in \sphm$. Since $\Phi$ is locally constant 
near the boundary and $\supp({\rm d}\Phi)$ is relatively compact in $M\backslash
C'$, there exist finitely many distinct points
$q_1,\cdots, q_k\in \sphm$ such that $\Phi(C')=\{q_1,\cdots,q_k \}$. Note that $(x,t)\in D_2$ implies $x\in C'$.
 We choose geodesic curves $\gamma_i:[0,1]\to \sphm$ such 
that $\gamma_i(0)=\star$ and $\gamma_i(1)=q_i$. 
Let $f:\mathbf{R}\to \mathbf{R}$ be a smooth cut-off function
such that if $t\leq (1-\varepsilon)/2$, $f(t)=1$ and if $t\geq (1+\varepsilon)/2$, $f(t)=0$.
%Let $p_1$, $p_2$ be the projection from $M\times [0,1]$ to 
%$M$, $[0,1]$ respectively.
Now we construct a map $\Theta: D(M)\to \sphm$ inspired by 
\cite{CZ24}*{Lemma~5.1} as follows\footnote{A similar trick also appears in \cite{LSWZ24+}.}:
\begin{equation}
\Theta\left((x,t)\right):=\begin{cases}
\Phi(x) \quad & \text{if}\  (x, t)\in M_1\times\{\frac{1}{2}-\varepsilon \},\\
\gamma_i(f(t)) \quad & \text{if} \  (x,t) \in D_2\ \text{and}\  x\in C\cap \Phi^{-1}(q_i), \\
\star \quad & \text{if}\  (x, t)\in M_1\times\{\frac{1}{2}+\varepsilon \}.
\end{cases}
\end{equation}
Observe that the restricted map $\Theta|_{D(M)\setminus (M_1\times\{\frac{1}{2}-\varepsilon \})}$ vanishes on $2$-vectors.
Then $\deg(\Theta)=\deg(\Phi)\neq 0$ and $\mathrm{supp}({\rm d}
\Theta)=\mathrm{supp}({\rm d}
\Phi)$.
%
%
%, we can construct a smooth map $\Psi: M^-\to N$ which is locally constant at infinity and of zero degree.
%Let $\Theta: {\rm d}M \to N$ be the smooth map defined by $\left.\Theta\right|_{M}=\Phi$ and $\left.\Theta\right|_{M^-}=\Psi$. 
Note that $\Theta$ is also locally constant at infinity. 

By contradiction, we assume that $\mathrm{scal}_{g_M}\geq 0$
on $M$. Since the scalar curvature on $D_2$ is positive
and the metric on $(x,t)\in D_1$ is the same as the metric on 
$x\in M$, we have $\mathrm{scal}_{g_{D(M)}}\geq 0$
and $\scal_{g_{D(M)}} \geq m(m-1)$  on $\supp({\rm d}\Theta)$,
which contradicts to \cref{thm:NNSC_sphere} by replacing $M$ by 
$D(M)$. This finishes the proof of \cref{thmx:complete_long_neck_NNSC}.
% \end{proof}

\begin{remark}
Consider the case where $(M,g_M)$ is only a noncompact
$m$-dimensional manifold with (possibly noncompact) boundary. 
Let $\Phi:M\to \sphm$ be a smooth area decreasing map which is locally constant near the boundary $\partial M$ and at infinity and of nonzero degree. 
Then $M\backslash \mathrm{supp}({\rm d}\Phi)$ has finitely many connected components, on which $\Phi$ is locally constant.
One of such connected components is called an \textit{end}.
Suppose that $\scal_{g_M} \geq m(m-1)$ on $\supp({\rm d}\Phi)$.

Suppose further that $M$ possesses at least one complete spin end without boundary and suppose $M$ admits the following
procedures of surgery to obtain a codimension zero complete noncompact spin submanifold $U$ with compact boundary $\partial U$:
\begin{enumerate}[label=\textup{(\roman*)}]
\item If an end is incomplete, we cut it off the manifold by a mean-convex
compact hypersurface such that the scalar curvature near which is positive.
\item If a remaining end has boundary, unless the boundary is compact, mean-convex and the scalar curvature near which is positive, we cut it off by a hypersurface
satisfying the same condition as above.
\item If the remaining part of the manifold is not spin, we cut the ends off
by hypersurfaces
satisfying the same condition as above until we obtain 
a complete noncompact spin manifold $U\subset M$ with compact boundary $\partial U$. 
See Figure \ref{figure}.
\end{enumerate} 
The submanifold $U$ obtained from the above surgery satisfies all assumptions in \cref{thmx:complete_long_neck_NNSC}.
Therefore, for such pair $(M, U)$ we have
\begin{equation}%\label{eq:NNSC_Bends}
	\inf_{x\in U} \scal_{g_M}(x)< 0.
\end{equation}

\tikzstyle{every node}=[scale=0.8]
\usetikzlibrary{patterns}
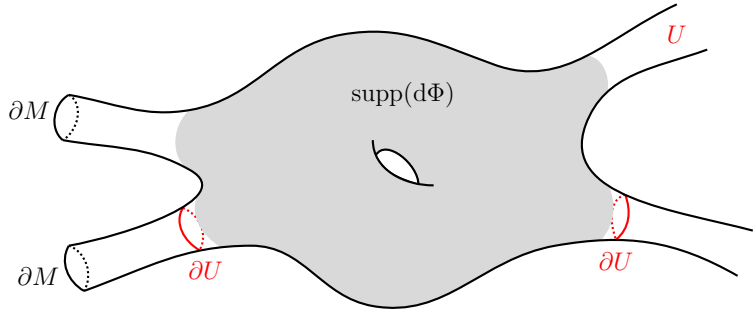
\begin{figure}[htpb]
	\centering
	\begin{tikzpicture}[scale=0.4]

	%%%%%%%%%%%%%%%%%%  left \partial U %%%%%%%%%%%%%%%%%
	
	% \draw[red,thick](4,3.5)to[out=240,in=150](4.5,2.3);
	% \draw[red,densely dotted,thick](4,3.5)to[out=-20,in=60](4.5,2.3);

	\draw[red,thick](5.75,4.25)to[out=240,in=150](6.3,2.85);
	\draw[red,densely dotted,thick](5.75,4.25)to[out=-20,in=60](6.3,2.85);
  
  %%%%%%%%%%%%%%%%%%  right \partial U %%%%%%%%%%%%%%%%%
	% \draw[red,thick](19.5,9.5)to[out=-5,in=80](20.2,8.5);
	% \draw[red,densely dotted,thick](19.5,9.5)to[out=250,in=190](20.2,8.5);
	
	% \draw[red,thick](22.5,4)to[out=-70,in=40](22.15,2.8); 
	% \draw[red,densely dotted,thick](22.5,4)to[out=200,in=110](22.15,2.8);

	\draw[red,thick](20.3,4.7)to[out=-70,in=40](20,3.2); 
	\draw[red,densely dotted,thick](20.3,4.7)to[out=200,in=110](20,3.2);

	% %%%%%%%%%%%%% support %%%%%%%%%%%%%%%%%%%%%%%%%
	\draw[fill=gray!30,gray!30](6,7.5)to[out=10,in=190](11,10)to[out=10,in=160](16,9)to[out=-20,in=210](19,9.3)to[out=-5,in=80](19.7,8.1)to[out=220,in=70](19,7)to[out=250,in=160](19.75,5)to[out=-70,in=50](19.5,3.2)to[out=185,in=10](18,3)to[out=190,in=-10](12,1)to[out=170,in=0](8,3)to[out=0,in=5](7,2.98)to[out=150,in=240](6.3,4.75)to[out=45,in=-30](6,5.5)to[out=150,in=220](6,7.5);
	
	%%%%%%%%%%%%%%% curves  %%%%%%%%%%%%%%%%%%%%%%%%%%%%%%%%%%%
	\draw[thick](2,8)to[out=-10,in=-170](6,7.5)to[out=10,in=190](11,10)to[out=10,in=160](16,9)to[out=-20,in=205](22,11);
	\draw[thick](1.8,6.5)to[out=-10,in=150](6,5.5)to[out=-30,in=10](2,3);
	\draw[thick](2.5,1.5)to[out=20,in=180](8,3)to[out=0,in=170](12,1)to[out=-10,in=190](18,3)to[out=10,in=150](24,2);
	\draw[thick](23,9.5)to[out=200,in=70](19,7)to[out=250,in=165](24.5,3.5);

	% %%%%%%%%%%%%%%%% compact boundary-1 %%%%%%%%%%%%%%%%%%%%%%%%%
	\draw[thick](2,8)to[out=210,in=140](1.8,6.5);
	\draw[densely dotted,thick](2,8)to[out=-30,in=30](1.8,6.5);
	
	% %%%%%%%%%%%%%%%% compact boundary-2 %%%%%%%%%%%%%%%%%%%%%%%%%
	\draw[thick](2,3)to[out=240,in=150](2.5,1.5);
	\draw[densely dotted,thick](2,3)to[out=-20,in=60](2.5,1.5);

	% %%%%%%%%%%%%%%%%%% colored hole %%%%%%%%%%%%%%%%%%%%%%%%%%%
	\draw[fill=white](12.1,6)to[out=60,in=90](13.5,5.03)to[out=160,in=-65](12.1,6);

	% %%%%%%%%%%%%%%%%%% hole %%%%%%%%%%%%%%%%%%%%%%%%%%%
	\draw[thick](12,6.5)to[out=-90,in=180](14,5);
	\draw[thick](12.1,6)to[out=60,in=90](13.5,5.03);

	%%%%%%%%%%%%%%%%%% right \Sigma %%%%%%%%%%%%%%%%%%%%%
	
%	\draw[yellow!40!green,thick](21,4.4)to[out=-70,in=40](20.6,3.15);
%	\draw[yellow!40!green,densely dotted,thick](21,4.4)to[out=200,in=110](20.6,3.15);

	%%%%%%%%%%%%%%%%%% left \Sigma %%%%%%%%%%%%%%%%%%%%%
	
%	\draw[yellow!40!green,thick](5.2,4)to[out=240,in=150](5.7,2.7);
%	\draw[yellow!40!green,densely dotted,thick](5.2,4)to[out=-20,in=60](5.7,2.7);

	%%%%%%%%%%%%%% labels %%%%%%%%%%%%%%%%%%%%%%%%%%%%%%%%%%%%%%%%%%%%%%%%%%%%%
	\node at (1,2) {$\partial M$};
	\node at (0.7,7.5) {$\partial M$};
	% \node at (9.5,10.5) {$\supp({\rm d}\Phi)$};
	% \node at (12,14) {$\textcolor{yellow!40!green}{M\setminus U_0}$};
	% \draw [->](12,13.5)--(12.5,10.5);
	\node at (22,10) {$\textcolor{red}{U}$};
	%\node at (20.5,2.5) {$\textcolor{yellow!40!green}{\Sigma}$};
	%\node at (6,2.2) {$\textcolor{yellow!40!green}{\Sigma}$};
	\node at (6.5,2.2) {$\textcolor{red}{\partial U}$};
	\node at (20,2.5) {$\textcolor{red}{\partial U}$};
	\node at (13,8) {${\rm supp}({\rm d}\Phi)$};
	\end{tikzpicture}
	\caption{An illustration of noncompact (necessarily not complete) manifold with surgery.}\label{figure}
\end{figure}
\end{remark}
    
%%%%%%%%%%%%%%%%%%%%%%%%%%%%%%%%%%%%%%%%%%%%%%%

\vspace{3mm}
\textbf{Acknowledgements.}
The authors are grateful to Professor Weiping Zhang and Professor Guangxiang Su for their helpful discussions and comments.
The first author
is partially supported by NSFC No.~12225105, and Science and Technology Commission of Shanghai Municipality (STCSM), grant No. 22DZ2229014.
The second author
is partially supported by
the National Natural Science Foundation of China grant No. 12501064, 
the China Postdoctoral Science Foundation (grant No. 2025M773075, Postdoctoral Fellowship Program grant No. GZC20252016 and Tianjin Joint Support Program grant No. 2025T002TJ)
and the Nankai Zhide Foundation.

%%%%%%%%%%%%%%%%%%%%%%%%%%%%%%%%%%%%%%%%%%%%%%%%%%%%%%%%%%%%%%%%%%%%%%%%%%%%%%
%%%%%%%%%%%%%%%%%%%%%%%%%%%%%%%%%%%%%%%%%%%%%%%%%%%%%%%%%%%%%%%%%%%%%%%%%%%%%%
%%%%%%%%%%%%%%%%%%%%%%%%%%%%%%%%%%%%%%%%%%%%%%%%%%%%%%%%%%%%%%%%%%%%%%%%%%%%%%
%%%%%%%%%%%%%%%%%%%%%%%%%%%%%%%%%%%%%%%%%%%%%%%%%%%%%%%%%%%%%%%%%%%%%%%%%%%%%%
%%%%%%%%%%%%%%%%%%%%%%%%%%%%%%%%%%%%%%%%%%%%%%%%%%%%%%%%%%%%%%%%%%%%%%%%%%%%%%
%%%%%%%%%%%%%%%%%%%%%%%%%%%%%%%%%%%%%%%%%%%%%%%%%%%%%%%%%%%%%%%%%%%%%%%%%%%%%%

% \bibliographystyle{amsalpha}


\begin{bibdiv}
\begin{biblist}


\bib{CZ24}{article}{
      author={Cecchini, S.},
      author={Zeidler, R.},
       title={Scalar and mean curvature comparison via the {D}irac operator},
        date={2024},
        ISSN={1465-3060},
     journal={Geom. Topol.},
      volume={28},
      number={3},
       pages={1167\ndash 1212},
         url={https://doi.org/10.2140/gt.2024.28.1167},
      review={\MR{4746412}},
}



\bib{deAl85}{article}{
   author={de Almeida, Sebasti\~ao},
   title={Minimal hypersurfaces of a positive scalar curvature manifold},
   journal={Math. Z.},
   volume={190},
   date={1985},
   number={1},
   pages={73--82},
   issn={0025-5874},
   review={\MR{0793350}},
   url={https://doi.org/10.1007/BF01159165},
}



% \bib{GL83}{article}{
%    author={Gromov, Mikhael},
%    author={Lawson, H. Blaine, Jr.},
%    title={Positive scalar curvature and the Dirac operator on complete
%    Riemannian manifolds},
%    journal={Inst. Hautes \'Etudes Sci. Publ. Math.},
%    number={58},
%    date={1983},
%    pages={83--196 (1984)},
%    issn={0073-8301},
%    review={\MR{0720933}},
% }



% \bib{Gro18}{article}{
%       author={Gromov, M.},
%        title={Metric inequalities with scalar curvature},
%         date={2018},
%         ISSN={1016-443X},
%      journal={Geom. Funct. Anal.},
%       volume={28},
%       number={3},
%        pages={645\ndash 726},
%       url={https://doi.org/10.1007/s00039-018-0453-z},
%       review={\MR{3816521}},
% }


\bib{Gr23}{article}{
   author={Gromov, Misha},
   title={Four lectures on scalar curvature},
   conference={
      title={Perspectives in scalar curvature. Vol. 1},
   },
   book={
      publisher={World Sci. Publ., Hackensack, NJ},
   },
   isbn={978-981-124-998-3},
   isbn={978-981-124-935-8},
   isbn={978-981-124-936-5},
   date={[2023] \copyright 2023},
   pages={1--514},
   review={\MR{4577903}},
}


\bib{GL80}{article}{
   author={Gromov, Mikhael},
   author={Lawson, H. Blaine, Jr.},
   title={Spin and scalar curvature in the presence of a fundamental group.
   I},
   journal={Ann. of Math. (2)},
   volume={111},
   date={1980},
   number={2},
   pages={209--230},
   issn={0003-486X},
   review={\MR{0569070}},
   url={https://doi.org/10.2307/1971198},
}


\bib{Ll98}{article}{
   author={Llarull, Marcelo},
   title={Sharp estimates and the Dirac operator},
   journal={Math. Ann.},
   volume={310},
   date={1998},
   number={1},
   pages={55--71},
   issn={0025-5831},
   review={\MR{1600027}},
   url={https://doi.org/10.1007/s002080050136},
}



\bib{LSWZ24+}{article}{
    author={Li, Y.},
    author={Su, G.},
    author={Wang, X.},
    author={Zhang, W.},
    title={Llarull's theorem on odd dimensional manifolds: the noncompact
    case},
    date={2024},
    eprint={https://arxiv.org/abs/2404.18153},
    url={https://arxiv.org/abs/2404.18153},
}



\bib{Shi25}{article}{
   author={Shi, Pengshuai},
   title={The odd-dimensional long neck problem via spectral flow},
   journal={Int. Math. Res. Not. IMRN},
   date={2025},
   number={17},
   pages={Paper No. rnaf262, 19},
   issn={1073-7928},
   review={\MR{4951381}},
   url={https://doi.org/10.1093/imrn/rnaf262},
}


% \bib{Shi25+}{article}{
%       author={Shi, P.},
%        title={Spectral flow of Callias operators, odd K-cowaist, and positive scalar curvature},
%         date={2025},
%       eprint={https://arxiv.org/abs/2501.17511},
%          url={https://arxiv.org/abs/2501.17511},
% }


\bib{Zh20}{article}{
   author={Zhang, Weiping},
   title={Nonnegative scalar curvature and area decreasing maps},
   journal={SIGMA Symmetry Integrability Geom. Methods Appl.},
   volume={16},
   date={2020},
   pages={Paper No. 033, 7},
   review={\MR{4089513}},
   url={https://doi.org/10.3842/SIGMA.2020.033},
}

\end{biblist}
\end{bibdiv}
\end{document}